\documentclass{gtart_h}  

\def\ifplaintex{\expandafter\ifx\csname documentclass\endcsname\relax}

\def\ifplaintex{\expandafter\ifx\csname documentclass\endcsname\relax}

\ifplaintex 
\else
\fi

\expandafter\ifx\csname epsfbox\endcsname\relax\input epsf\fi

\def\gt{{\mathsurround=0pt\it $\cal G\mskip-2mu$eometry \&\ 
$\cal T\!\!$opology}}        

\def\gtp{{\mathsurround=0pt\it $\cal G\mskip-2mu$eometry \&\ 
$\cal T\!\!$opology $\cal P\!$ublications}}  

\def\lognumber#1{\def\thelognumber{#1}}
\def\volumenumber#1{\def\thevolumenumber{#1}}
\def\papernumber#1{\def\thepapernumber{#1}}
\def\volumeyear#1{\def\thevolumeyear{#1}}

\def\pagenumbers#1#2{\def\startpage{#1}\def\finishpage{#2}}
\def\published#1{\def\publishdate{#1}}
\def\proposed#1{\def\theproposer{#1}}
\def\seconded#1{\def\theseconders{#1}}
\def\received#1{\def\receiveddate{#1}}
\def\revised#1{\def\reviseddate{#1}}
\def\accepted#1{\def\accepteddate{#1}}

\def\asciiaddress#1{\def\theasciiaddress{#1}}
\def\asciiemail#1{\def\theasciiemail{#1}}

\long\def\asciiabstract#1{\long\def\theasciiabstract{#1}}
\def\asciikeywords#1{\def\theasciikeywords{#1}}

\let\\\par\let\thelognumber\relax
\let\thevolumenumber\relax\let\thepapernumber\relax
\let\thevolumeyear\relax\let\thesamplenumber\relax\let\startpage\relax
\let\finishpage\relax\let\publishdate\relax\let\receiveddate\relax
\let\reviseddate\relax\let\accepteddate\relax\let\theasciititle\relax
\let\theasciiauthors\relax\let\theasciiaddress\relax
\let\theasciiabstract\relax\let\theasciikeywords\relax
\let\theasciiemail\relax\let\theshortauthors\relax\let\theshorttitle\relax

\long\def\maketitlep{   

\count0=\startpage

\gt\hfill      
\hbox to 77pt{\vbox to 0pt{\vglue -15pt\epsfbox{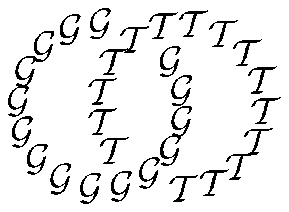}\vss}\hss}
\break
{\small\ifx\thesamplenumber\relax 
Volume \else Sample
\fi\thevolumenumber\ (\thevolumeyear)
\startpage--\finishpage\nl
Published: \publishdate}
\vglue 0.5truein plus 0.4fil minus 0.1truein

{\parskip=0pt\leftskip 0pt plus 1fil\def\\{\par\smallskip}{\ifplaintex\large
\else\Large\fi\bf\thetitle}\par\medskip}   

\vglue 0pt plus 0.1fil 

{\parskip=0pt\leftskip 0pt plus 1fil\def\\{\par}{\sc\theauthors}
\par\medskip}

\vglue 0pt plus 0.1fil 

{\small\parskip=0pt\let\newline\\
{\leftskip 0pt plus 1fil\def\\{\par}{\sl\theaddress}\par}
\expandafter\ifx\theemail\relax    
\relax\else\vglue 5pt plus 0.02fil minus 2pt\def\\{\stdspace{\rm 
and}\stdspace} 
\cl{Email:\stdspace\tt\theemail}\fi
\ifx\theurl\relax                  
\relax\else\vglue 5pt plus 0.02fil minus 2pt\def\\{\stdspace{\rm 
and}\stdspace}
\cl{URL:\stdspace\tt\theurl}\fi\par}

\vglue 7pt plus 0.3fil minus 3pt

{\bf Abstract}
\vglue 5pt plus 0.1fil minus 2pt

\theabstract

\vglue 7pt plus 0.3fil minus 3pt

{\bf AMS Classification numbers}\quad Primary:\quad \theprimaryclass

Secondary:\quad \thesecondaryclass

\vglue 5pt plus 0.3fil minus 2pt

{\bf Keywords:}\quad \thekeywords

\vglue 10pt plus 0.5fil minus 5pt

{\small  Proposed: \theproposer\hfill Received: \receiveddate\nl
Seconded: \theseconders\hfill 
\ifx\reviseddate\relax                         
Accepted: \accepteddate                        
\else
Revised: \reviseddate                          
\fi}
\eject
}       

\font\phead=cmsl9 scaled 950
\font=cmsl9 scaled 1050
\font\pnum=cmbx10 scaled 913
\font\lnum=cmbx10 
\font\pfoot=cmsl9 scaled 950
\font=cmsl9 scaled 1050
\ifplaintex
\headline{\vbox to 0pt{\vskip -4.5mm\line{\small\phead\ifnum
\count0=\startpage ISSN 1364-0380 (on line)
1465-3060 (printed) \hfill {\pnum\folio}\else\ifodd\count0\def\\{ }%
\ifx\theshorttitle\relax\thetitle\else\theshorttitle\fi\hfill{\pnum\folio}
\else\def\\{ and }{\pnum\folio}\hfill\ifx\theshortauthors\relax\theauthors
\else\theshortauthors\fi\fi\fi}\vss}}
\footline{\vbox to 0pt{\vglue 0mm\line{\small\pfoot\ifnum\count0=\startpage
\copyright\ \gtp\hfill\else
\gt, Volume \thevolumenumber\ (\thevolumeyear)\hfill\fi}\vss
}}
\else
\makeatletter
\def\@oddhead{{\small\lhead\ifnum\count0=\startpage ISSN 1364-0380 (on line)
1465-3060 (printed) \hfill {\lnum\number\count0}\else\ifodd\count0
\def\\{ }\ifx\theshorttitle\relax \thetitle \else\theshorttitle\fi\hfill
{\lnum\number\count0}\else\def\\{ and }{\lnum\number\count0}
\hfill\ifx\theshortauthors\relax 
\theauthors\else\theshortauthors\fi\fi\fi}}\def\@evenhead{\@oddhead}
\def\@oddfoot{\small\lfoot\ifnum\count0=\startpage\copyright\ \gtp\hfill\else
\gt, Volume \thevolumenumber\ (\thevolumeyear)\hfill\fi}
\def\@evenfoot{\@oddfoot}
\makeatother
\fi

\newwrite\gtoutfile
\long\gdef\makeheadfile{  
{\def\\{, }\def\s{ }
\immediate\openout\gtoutfile head.xxx
\immediate\write\gtoutfile{Proxy-for: \ifx\theasciiauthors\relax
\theauthors\else\theasciiauthors\fi\s<\ifx\theasciiemail\relax\theemail\else\theasciiemail\fi>}
\immediate\write\gtoutfile{\noexpand\\}
\immediate\write\gtoutfile{Authors: \ifx\theasciiauthors\relax
\theauthors\else\theasciiauthors\fi}
{\def\\{ }\immediate\write\gtoutfile{Title: \ifx\theasciititle\relax
\thetitle\else\theasciititle\fi}}
\immediate\write\gtoutfile{Subj-class: GT or SG or MG etc}
\immediate\write\gtoutfile{MSC-class: \theprimaryclass\ifx\thesecondaryclass\relax\else, \thesecondaryclass\fi}
\immediate\write\gtoutfile{Journal-ref: Geom. Topol. \thevolumenumber
(\thevolumeyear) \startpage-\finishpage}
\immediate\write\gtoutfile{Comments: Published by Geometry and Topology at}
\immediate\write\gtoutfile{\s\s http://www.maths.warwick.ac.uk/gt/GTVol\thevolumenumber/paper\thepapernumber.abs.html}
\immediate\write\gtoutfile{\noexpand\\}
\immediate\write\gtoutfile{}
\ifx\theasciiabstract\relax
\immediate\write\gtoutfile{\theabstract}\else
\immediate\write\gtoutfile{\theasciiabstract}\fi
\immediate\write\gtoutfile{}
\immediate\write\gtoutfile{\noexpand\\}
\immediate\write\gtoutfile{}
\immediate\closeout\gtoutfile}}  

\def\maketitlepage{\maketitlep\makeheadfile}
\let\maketitle\maketitlepage

\lognumber{359}
\volumenumber{8}\papernumber{25}\volumeyear{2004}
\pagenumbers{947}{968}
\received{4 September 2003}
\published{29 June 2004}
\revised{19 April 2004}
\accepted{3 June 2004}
\proposed{Peter Kronheimer}
\seconded{Robion Kirby, Yasha Eliashberg}

\usepackage{amsmath,amssymb,amscd,curves}

\theoremstyle{plain}
\newtheorem{thm}{Theorem}[section]
\newtheorem{lem}[thm]{Lemma}
\newtheorem{cor}[thm]{Corollary}
\newtheorem{prop}[thm]{Proposition}

\theoremstyle{definition}
\newtheorem{defn}[thm]{Definition}

\def \S {Section }

\def \R {\mathbf{R}}
\def \Z {\mathbf{Z}}
\def \Sig{\Sigma}

\def \vp {\varphi}

\def \a {\alpha}
\def \b {\beta}
\def \g {\gamma}
\def \d {\delta}
\def \k {\kappa}
\def \lam {\lambda}

\def \o {\omega}
\def \s {\sigma}
\def \t {\tau}

\def \bd {\partial}
\def \x {\times}
\def \- {\setminus}
\def \C {\subset}
\def \ve {\varepsilon}

\def\pt{\text{pt}}
\def \bX {\bar{X}}

\def \ssw {\text{SW}}
\def \sw {\mathcal{SW}}
\def \swb{\overline {\mathcal{SW}}}
\def \swu{\text{SW}'}

\def \DD {\Delta}

\def \Ds {\DD^{\text{\it{sym}}}}

\def\I{\mathcal{I}}
\def\tm{T_{\mu}}

\def\fs{\mathfrak{s}}
\def\Ci{C^{\infty}}

\begin{document}

\title{Invariants for Lagrangian tori}
\authors{Ronald Fintushel\\Ronald J Stern}
\address{Department of Mathematics, Michigan State University\\
East Lansing, Michigan 48824, USA}
\secondaddress{Department of Mathematics, University of California \\
Irvine, California 92697, USA}

\asciiaddress{Department of Mathematics, Michigan State
University\\East Lansing, Michigan 48824, USA\\and\\Department of
Mathematics, University of California\\Irvine, California 92697, USA}

\gtemail{\mailto{ronfint@math.msu.edu}{\rm\qua 
and\qua}\mailto{rstern@math.uci.edu}}
\asciiemail{ronfint@math.msu.edu, rstern@math.uci.edu}

\begin{abstract}  We define an simple invariant $\lam(T)$ of an embedded nullhomologous Lagrangian torus and use this invariant to show that many symplectic $4$--manifolds have infinitely many pairwise symplectically inequivalent nullhomologous Lagrangian tori. We further show that for a large class of examples that $\lam(T)$ is actually a $\Ci$ invariant. In addition, this invariant is used to show that many symplectic $4$--manifolds have nontrivial homology classes which are represented by infinitely many pairwise inequivalent Lagrangian tori, a result first proved by S Vidussi for the homotopy K3--surface obtained from knot surgery using the trefoil knot \cite{V}.
\end{abstract}

\asciiabstract{%
We define an simple invariant of an embedded nullhomologous Lagrangian
torus and use this invariant to show that many symplectic 4-manifolds
have infinitely many pairwise symplectically inequivalent
nullhomologous Lagrangian tori. We further show that for a large class
of examples that lambda(T) is actually a C-infinity invariant. In
addition, this invariant is used to show that many symplectic
4-manifolds have nontrivial homology classes which are represented by
infinitely many pairwise inequivalent Lagrangian tori, a result first
proved by S Vidussi for the homotopy K3-surface obtained from knot
surgery using the trefoil knot in [Lagrangian surfaces in a fixed
homology class: existence of knotted Lagrangian tori, J. Diff. Geom. 
(to appear)].}

\primaryclass{57R57}
\secondaryclass{57R17}

\keywords{$4$--manifold, Seiberg--Witten invariant, symplectic, Lagrangian}
\asciikeywords{4-manifold, Seiberg-Witten invariant, symplectic, Lagrangian}

\maketitlepage

\section{Introduction\label{Intro}}

This paper is concerned with the construction and detection of homologous but inequivalent (under diffeomorphism or symplectomorphism) Lagrangian tori. In recent years there has been considerable progress made on the companion problem for symplectic tori. In fact, it is now known that in a simply connected symplectic $4$--manifold with $b^+>1$, if $T$ is an embedded symplectic torus, then for each $m>1$ the homology class of $mT$ contains infinitely many nonisotopic embedded symplectic tori. (See eg \cite{EP, surfaces, tori}.) The techniques used for the construction of these examples fail for Lagrangian tori.

The first examples of inequivalent homologous Lagrangian tori were discovered by  S Vidussi \cite{V} who presented a technique for constructing infinitely many homologous, but nonisotopic, Lagrangian tori in $E(2)_K$, the result of knot surgery on the K3--surface using the trefoil knot, $K$. The work in this paper was motivated by an attempt to better understand and distinguish the examples presented in \cite{V}. We found that the key to these examples is the construction of infinite families of nullhomologous Lagrangian tori $T$ in a symplectic $4$--manifold $X$. There is a simple process by which an integer $\lam(T)$, a Lagrangian framing defect, can be associated to $T$. In this paper we show that $\lam(T)$ is an invariant of the symplectomorphism, and in many cases the diffeomorphism,  type of $(X,T)$. We then construct  infinite families of inequivalent  nullhomologous Lagrangian tori distinguished by $\lam(T)$.  Homologically essential examples are created from these by a circle sum process. 

Some examples typical of those we which we study can be briefly described: Let $X$ be any symplectic manifold which contains an embedded self-intersection $0$ symplectic torus $T$. For any fibered knot $K$ consider the symplectic manifold $X_K$ constructed by knot surgery \cite{KL4M}. Since  $X_K$ is the fiber sum  of $X$  and $S^1 \x M_K$ along $T$ and $S^1\times m$ in $S^1 \x M_K$,  where $M_K$ is the result of $0$--framed surgery on $K$ and $m$ is a meridian of $K$, there is a codimension~$0$ submanifold $V$ in $X_K$ diffeomorphic to $S^1 \x (M_K\- m)$. The manifold $V$ is fibered by punctured surfaces $\Sig'$, and if $\g$ is any loop on such a surface, the torus $T_\g=S^1\x\g$ in $X_K$ is nullhomologous and Lagrangian.  We show that $\lam(T_\g)$ is a diffeomorphism invariant. This invariant persists even after circle sums with essential Lagrangian tori, and it distinguishes all our (and Vidussi's) examples.

Here is a more precise summary of our examples:

\begin{thm} \label{T1}{\rm{(a)}}\qua Let $X$ be any symplectic manifold with $b_2^+(X)>1$ which contains an embedded self-intersection $0$ symplectic torus with a vanishing cycle. (See \S\ref{NHT} for a definition.) Then for each nontrivial fibered knot $K$ in $S^3$, the result of knot surgery $X_K$ contains infinitely many nullhomologous Lagrangian tori, pairwise inequivalent under orientation-preserving diffeomorphisms.

\noindent{\rm{(b)}}\qua Let $X_i$, $i=1,2$, be symplectic $4$--manifolds containing embedded self-intersection $0$ symplectic tori $F_i$ and assume that $F_1$ contains a vanishing cycle. Let $X$ be the fiber sum, $X=X_1\#_{F_1=F_2}X_2$. Then for each nontrivial fibered knot $K$ in $S^3$, the manifold $X_K$ contains an infinite family of homologically primitive and homologous Lagrangian tori which are pairwise inequivalent.
\end{thm}

In Section \ref{SELT} we shall also give examples of nullhomologous Lagrangian tori $T_i$ in a symplectic $4$--manifold where the $\lam(T_i)$ are mutually distinct, so these tori are inequivalent under symplectomorphisms, but the techniques of this paper, namely relative Seiberg--Witten invariants, fail to distinguish the the $T_i$. It remains an extremely interesting question whether the these tori are equivalent under diffeomorphisms.

\section{Seiberg--Witten invariants for embedded tori}

The Seiberg--Witten
invariant  of a smooth closed oriented $4$--manifold
$X$ with $b_2^+(X)>1$ is an integer-valued function $\ssw_X$ which is defined on the set of $spin ^c$ structures over $X$. Corresponding to each $spin ^c$ structure $\fs $ over $X$ is the bundle of positive spinors $W^+_{\fs}$ over $X$. Set $c(\fs)\in H_2(X)$ to be the Poincar\'e dual of $c_1(W^+_{\fs})$.
Each $c(\fs)$ is a characteristic element of $H_2(X;\Z)$ (ie, its
Poincar\'e dual $\hat{c}(\fs)=c_1(W^+_{\fs})$ reduces mod~2 to $w_2(X)$).
We shall work with the modified
Seiberg--Witten invariant 
\[ \swu_X\co  \lbrace k\in H_2(X;\Z)|\hat{k}\equiv w_2(TX)\pmod2)\rbrace
\rightarrow \Z \]
defined by $\swu_X(k)=\sum\limits_{c(\fs)=k}\ssw_X(\fs)$.

The sign of $\ssw_X$
depends on a homology orientation of $X$, that is, an orientation of
$H^0(X;\R)\otimes\det H_+^2(X;\R)\otimes \det H^1(X;\R)$. If $\swu_X(\b)\neq
0$, then $\b$ is called a {\it{basic class}} of $X$. It is a fundamental
fact that the set of
basic classes is
finite. Furthermore, if $\b$ is a basic class, then so is $-\b$ with
$\swu_X(-\b)=(-1)^{(\text{e}+\text{sign})(X)/4}\,\swu_X(\b)$ where
$\text{e}(X)$ is
the Euler number and $\text{sign}(X)$ is the signature of $X$. The Seiberg--Witten invariant is an orientation-preserving diffeomorphism invariant of $X$ (together with the choice of a homology orientation).

It is convenient to view the Seiberg--Witten invariant as an element of the integral group ring $\Z H_2(X)$, where for each $\a\in H_2(X)$ we let $t_\a$ denote the corresponding element in $\Z H_2(X)$. Suppose that 
$\{\pm \b_1,\dots,\pm \b_n\}$ is the set of nonzero basic classes for $X$.
Then the Seiberg--Witten invariant of $X$ is the Laurent polynomial
\[\sw_X = \swu_X(0)+\sum_{j=1}^n \swu_X(\b_j)\cdot
(t_{\b_j}+(-1)^{(\text{e}+\text{sign})(X)/4}\, t_{\b_j}^{-1}) \in \Z H_2(X).\]
Suppose that $T$ is an embedded (but not necessarily homologically essential) torus of self-intersection $0$ in $X$, and identify a tubular neighborhood of $T$ with $T\x D^2$. 
Let $\a$, $\b$, $\g$ be simple loops on $\bd(T\x D^2)$ whose homology classes generate $H_1(\bd(T\x D^2))$. Denote by $X_T(p,q,r)$ the result of surgery on $T$ which annihilates the class of $p\a +q\b + r\g$; ie,
\begin{equation}\label{surgery}
X_T(p,q,r) = (X\- T\x D^2) \cup_{\vp} T^2\x D^2
\end{equation}
where $\vp\co \bd(X\- T\x D^2)\to \bd(T^2\x D^2)$ is an orientation-reversing diffeomorphism satisfying $\vp_*[p\a +q\b + r\g] =[\bd D^2]$. An important formula for calculating the Seiberg--Witten invariants of surgeries on tori is due to Morgan, Mrowka, and Szabo \cite{MMS} (see also \cite{MT}, \cite{T}). Suppose that $b_2^+(X\- (T\x D^2))>1$. Then each $b_2^+(X_T(p,q,r))>1$. Given a class $k\in H_2(X)$:
\begin{multline}\label{surgery formula} \sum_i\swu_{X_T(p,q,r)}(k_{(p,q,r)}+i[T])= 
p\sum_i\swu_{X_T(1,0,0)}(k_{(1,0,0)}+i[T]) +\\+q\sum_i\swu_{X_T(0,1,0)}(k_{(0,1,0)}+i[T])
+r\sum_i\swu_{X_T(0,0,1)}(k_{(0,0,1)}+i[T])
\end{multline}
In this formula, $T$ denotes the torus which is the core $T^2\x {0}\C T^2\x D^2$ in each specific manifold $X_T(a,b,c)$ in the formula, and $k_{(a,b,c)}\in H_2(X_T(a,b,c))$ is any class which agrees with the restriction of $k$ in $H_2(X\- T\x D^2,\bd)$ in the diagram:
\[ \begin{array}{ccc}
H_2(X_T(a,b,c)) &\longrightarrow & H_2(X_T(a,b,c), T\x D^2)\\
&&\Big\downarrow \cong\\
&&H_2(X\- T\x D^2,\bd)\\
&&\Big\uparrow \cong\\
H_2(X)&\longrightarrow & H_2(X,T\x D^2)
\end{array}\]
Furthermore, in each term of \eqref{surgery formula}, unless the homology class $[T]$ is 2--divisible, each $i$ must be even since the classes $k_{(a,b,c)}+i[T]$ must be characteristic in $H_2(X_T(a,b,c))$.

Let $\pi(a,b,c)\co  H_2(X_T(a,b,c))\to H_2(X\- T\x D^2,\bd)$ be the composition of maps in the above diagram, and $\pi(a,b,c)_*$ the induced map of integral group rings. Since we are interested in invariants of the pair $(X,T)$, we shall work with \[\swb_{(X_T(a,b,c),T)}=\pi(a,b,c)_*(\sw_{X_T(a,b,c)})\in \Z H_2(X\- T\x D^2,\bd).\]
The indeterminacy in \eqref{surgery formula} is caused by multiples of $[T]$; so
passing to $\swb$ removes this indeterminacy, and the Morgan--Mrowka--Szabo formula becomes
\begin{equation}\label{surgery formula 2}
 \swb_{(X_T(p,q,r),T)} = 
p\,\swb_{(X_T(1,0,0),T)} +q\,\swb_{(X_T(0,1,0),T)}
+r\,\swb_{(X_T(0,0,1),T)}.
\end{equation}
Let $T$ and $T'$ be embedded tori in the oriented $4$--manifold $X$. We shall say that these tori are {\it $\Ci$\!--equivalent} if there is an orientation-preserving diffeomorphism  $f$ of $X$ with $f(T)=T'$. Any self-diffeomorphism of $X$ which throws $T$ onto $T'$, takes a loop on the $\bd(T\x D^2)$ to a loop on the boundary of a tubular neighborhood of $T'$. Set  \[\I(X,T) = \{\swb_{(X_T(a,b,c),T)}|a,b,c\in\Z\}.\]

\begin{prop} Let $T$ be an embedded torus of self-intersection $0$ in the simply connected $4$--manifold $X$ with $b_2^+(X\- T)>1$. After fixing a homology orientation for $X$, $\I(X,T)$ is an invariant of the pair $(X,T)$ up to $\Ci$\!-equivalence.\qed
\end{prop}

\section{The Lagrangian framing invariant}

In this section we shall define the invariant $\lam(T)$ of a nullhomologous Lagrangian torus. To begin, consider a nullhomologous torus $T$ embedded in a smooth $4$--manifold $X$ with tubular neighborhood $N_T$. Let $i\co \bd N_T\to X\- N_T$ be the inclusion. 

\begin{defn} A {\it{ framing}} of $T$ is a diffeomorphism $\vp\co  T\x D^2\to N_T$ such that $\vp(p)=p$ for all $p\in T$. A framing $\vp$ of $T$ is {\it{nullhomologous}} if for $x\in \bd D^2$, the homology class 
$\vp_*[T\x \{x\}]\in \ker i_*$.
\end{defn}

Given a framing $\vp\co  T\x D^2\to N_T$, there is an associated section $\s(\vp)$ of $\bd N_T\to T$ given by $\s(\vp)(x)=\vp(x,1)$, and given a pair of framings, $\vp_0$, $\vp_1$ there is a difference class $\d(\vp_0,\vp_1)\in H^1(T;\Z)\cong [T, S^1]$, the homotopy class of the composition 
\[ T \stackrel{\s(\vp_1)} {\longrightarrow} \bd N_T \stackrel{\vp_0^{-1}} {\longrightarrow} T\x \bd D^2 \stackrel{\text{pr}_2} {\longrightarrow} \bd D^2 \cong S^1. \]
Note that if $\a$ is a loop on $T$, then $\d(\vp_0,\vp_1)[\a]=
\s(\vp_1)_*[\a]\cdot\s(\vp_0)_*[T]$ using the intersection pairing on $\bd N_T$, or equivalently, 
\[ (\vp_0)_*^{-1}\s(\vp_1)_*[\a] = [\a\x \{1\}] + \d(\vp_0,\vp_1)[\a]\, [\bd D^2] \in H_1(T\x \bd D^2;\Z). \]

\begin{prop} A nullhomologous framing of $T$ is unique up to homotopy.
\end{prop}
\begin{proof}  Since $T$ is nullhomologous, it follows that $H_2(X\- N_T)\to H_2(X)$
is onto, and $H_3(X\- N_T,\bd N_T)\cong H_3(X,T)\cong H_3(X)\oplus H_2(T)$. Then the long exact sequence of $(X\- N_T, \bd N_T)$ shows that the kernel of 
$i_*\co H_2(\bd N_T)\to H_2(X\- N_T)$ is isomorphic to $H_2(T)=\Z$. 
So any two nullhomologous framings $\vp_0$, $\vp_1$ give rise to homologous tori $\s(\vp_i)(T)$ in $\bd N_T$. Thus for any loop $\a$ on $T$: 
\begin{multline*} \d(\vp_0,\vp_1)[\a]=\s(\vp_1)_*[\a]\cdot\s(\vp_0)_*[T] = \\ \s(\vp_1)_*[\a]\cdot\s(\vp_1)_*[T] = [\a \x \{1 \}]\cdot [T\x \{ 1\}] = 0
\end{multline*}
the last pairing in $T\x \bd D^2$. Hence $\d(\vp_0,\vp_1)=0$.
\end{proof}

\noindent We denote by $\vp_N$ any such nullhomologous framing of $T$.

Now suppose that $(X,\o)$ is a symplectic $4$--manifold containing an embedded Lagrangian torus $T$.  For any closed oriented Lagrangian surface $\Sig\C X$ there is a nondegenerate bilinear pairing 
\[ (TX/T\Sig) \otimes T\Sig \to \R, \ \ \ \ ([v],u)\to \o(v,u).\] 
Hence, the normal bundle $N_{\Sig}\cong T^*\Sig$, the cotangent bundle; so computing Euler numbers, 
$ 2g-2 = -e(\Sig)= e(T^*(\Sig))=e(N_{\Sig})=\Sig\cdot\Sig$, for $g$ the genus of $\Sig$. Furthermore, this is true symplectically as well. The Lagrangian neighborhood theorem \cite{W} states that each such Lagrangian surface has a tubular neighborhood which is symplectomorphic to a neighborhood of the zero section of its cotangent bundle with its standard symplectic structure, where the symplectomorphism is the identity on $\Sig$. 

Thus an embedded Lagrangian torus $T$ has self-intersection $0$, and small enough tubular neighborhoods $N_T$ have, up to symplectic isotopy, a preferred framing  
$\vp_L\co T\x D^2 \to N_T$
such that for any point $x \in D^2$, the torus $\vp_L(T\x \{ x\})$ is also Lagrangian. We shall call $\vp_L$ the {\it Lagrangian framing} of $T$. 

Thus if $T$ is a nullhomologous Lagrangian torus, we may consider the difference
$\d(\vp_N,\vp_L)\in H^1(T;\Z)=[T,S^1]$. It thus induces a well-defined homomorphism $\d(\vp_N,\vp_L)_*\co  H_1(T;\Z)\to H_1(S^1)$.

\begin{defn} The {\it{Lagrangian framing invariant}} of a nullhomologous Lagrangian torus $T$ is the nonnegative integer $\lam(T)$ such that 
$$\d(\vp_N,\vp_L)_*(H^1(T;\Z))=\lam(T)\, \Z {\rm\qua in\qua} H_1(S^1)=\Z.$$
\end{defn}
\noindent Thus if $a$ and $b$ form a basis for $H_1(T;\Z)$ then $\lam(T)$ is the greatest common divisor of $|\d(\vp_N,\vp_L)(a)|$ and $|\d(\vp_N,\vp_L)(b)|$. Furthermore, if $f\co X\to Y$ is a symplectomorphism with $f(T) = T'$, then $f\circ\vp_N$ is a nullhomologous framing of $T'$, and for small enough $N_T$, $f\circ\vp_L$ is the Lagrangian framing of $T'$. Hence:

\begin{thm} Let $T$ be a nullhomologous Lagrangian torus in the symplectic $4$--manifold $X$. Then the Lagrangian framing invariant $\lam(T)$ is a symplectomorphism invariant of $(X,T)$. \qed
\end{thm}

Here is an example. Let $K$ be any fibered knot in $S^3$, and let $M_K$ be the result of $0$--surgery on $K$. Then $M_K$ is a $3$--manifold with the same homology as $S^2\x S^1$,  and $M_K$ is fibered over the circle. Let $\g$ be any  embedded loop which lies on a fiber of the fibration $S^3\- K\to S^1$. Ie, $\g$ lies on a Seifert surface $\Sig$ of $K$. The first homology $H_1(M_K)=H_1(S^3\- K)\cong \Z$, and the integer corresponding to a given loop is the linking number of the loop with $K$. Since $\g$ lies on a Seifert surface, its linking number with $K$ is $0$, and so $\g$ is nullhomologous in $M_K$.

Taking the product with a circle, $S^1\x M_K$ fibers over $T^2$, and it is a symplectic $4$--manifold with a symplectic form which arises from the sum of volume forms in the base and in the fiber. More precisely, one can choose metrics so that the fiber bundle projection, $p\co M_K\to S^1$ is harmonic. Let $\a$ be the volume form on the base $S^1$, and let $\b$ be the volume form on the first $S^1$ in $S^1\x M_K$. Then $\o=\b\wedge p^*(\a)+*_3\,p^*(\a)$ defines a symplectic form on $S^1\x M_K$. Since $\g$ lies in a fiber, its tangent space at any point is spanned by a vector parallel to the tangent space of the fiber and a vector tangent to $S^1$. So $\o$ vanishes on $T(S^1\x \g)$. (See also \cite{V}.) 
Thus $T_\g = S^1\x \g$ is a nullhomologous Lagrangian torus in $S^1\x M_K$.
If $\g'$ is a pushoff of $\g$ in the Seifert surface $\Sig$, then $S^1\x \g'$ is again Lagrangian. This, together with the pushoff of $\g$ onto nearby fibers, describes the Lagrangian framing of $T_\g$. We shall also say that $\g'$ is the {\it Lagrangian pushoff} of $\g$.

\begin{defn} Let $K$ be a fibered knot in $S^3$, and let $\g$ be any embedded loop lying on a fiber of the fibration $S^3\- K\to S^1$. The {\it Lagrangian framing defect} $\lam(\g)$ of $\g$ is the linking number of $\g$ with a  Lagrangian pushoff of itself.
\end{defn}
\noindent In Figures 1 and 2, we have $\lam(\g_1)=1$, and $\lam(\g_2)=3$. 

\centerline{\small\unitlength 1cm
\begin{picture}(12,6)
\put (2,3){\oval(3,4)[l]}
\put (3,3){\oval(3,4)[r]}
\put (2,5){\line(3,-4){1}}
\put (3,3.66){\line(-3,-4){.4}}
\put (3,5){\line(-3,-4){.4}}
\put (2,3.66){\line(3,4){.4}}
\put (2,3.66){\line(3,-4){1}}
\put (2,2.33){\line(3,4){.4}}
\put (2,2.33){\line(3,-4){1}}
\put (3,2.33){\line(-3,-4){.4}}
\put (2,1){\line(3,4){.4}}
\put (2.4,3.66){\oval(1.5,1.4)[l]}
\put (2.6,3.66){\oval(1.5,1.4)[r]}
\put (1.5,4.25){\small{$\g_1$}}
\put (1.85,.5){Figure 1}
\put (9,3){\oval(3,4)[l]}
\put (10,3){\oval(3,4)[r]}
\put (9,5){\line(3,-4){1}}
\put (10,3.66){\line(-3,-4){.4}}
\put (10,5){\line(-3,-4){.4}}
\put (9,3.66){\line(3,4){.4}}
\put (9,3.66){\line(3,-4){1}}
\put (9,2.33){\line(3,4){.4}}
\put (9,2.33){\line(3,-4){1}}
\put (10,2.33){\line(-3,-4){.4}}
\put (9,1){\line(3,4){.4}}
\put (9.4,3.69){\oval(1.5,1.3)[l]}
\put (9.6,2.33){\oval(1.45,1.5)[r]}
\qbezier(9.58,2.98)(9.75,3)(9.8,3.02)
\qbezier(9.95,3.1)(11.1,3.8)(9.58,4.35)
\put (9.48,3.07){\line(1,0){.15}}
\put (9.4,2.26){\oval(1.5,1.4)[l]}
\put (8.5,4.25){\small{$\g_2$}}
\put (8.85,.5){Figure 2}
\put(9.38,1.56){\line(1,0){.125}}
\put(9.38,2.96){\line(1,0){.1}}
\end{picture}}

As in the above definition, let $\g$ be an embedded loop lying on a fiber of the fibration $S^3\- K\to S^1$, and let $N(\g)\cong \g\x D^2$ be a tubular neighborhood of $\g$. Further, let $\ell(\g)\in H_1(\bd N(\g))$ be the (nullhomologous) $0$--framing of $\g$ in $S^3$; that is, the nontrivial primitive class which is sent to $0$ by $H_1(\bd(N(\g)))\to H_1(S^3\- N(\g))$. Then if $\g'$ is a Lagrangian pushoff of $\g$, in $H_1(\bd(N(\g)))$ we have the relation 
\[ [\g'] = \ell(\g) + \lam(\g) [\bd D^2]. \]
In other words, the Lagrangian pushoff corresponds to the framing $\lam(\g)$ with respect to the usual $0$--framing $\ell(\g)$ in $S^3$. So, for example, a Lagrangian  $1/p$  surgery on  the curve $\g_1$ above corresponds to a $(p+1)/p$ surgery with respect to the usual framing of $\g_1$ in $S^3$. More generally, a $1/p$ Lagrangian surgery on a curve $\g$ in the Seifert surface of a fibered knot in $S^3$ corresponds to a $(p\lam(\g)+1)/p$ surgery with respect to the usual framing of $\g$ in $S^3$.

\begin{thm} \label{LFI} In $S^1\x M_K$, the Lagrangian framing invariant of $T_\g$ is $\lam(T_\g) =|\lam(\g)|$.
\end{thm}
\begin{proof} As a basis for $H_1(T_\g;\Z)$ take $[\{1\}\x\g]$ and 
$[S^1\x\{ x\}]$ where $x\in \g$. Since the linking number of $\g$ and $K$ is $0$, 
there is a Seifert surface $C$ for $\g$ in $S^3$ which is disjoint from $K$. The tubular neighborhood $N_T$ of $T_\g$ is given by $N_T=S^1\x N(\g)$ and 
$\s(\vp_N)(T_\g)=(S^1\x C)\cap \bd N_T = S^1\x (C\cap\bd N(\g))= S^1\x \ell(\g)$,
where we are using `$\ell(\g)$' here to denote a curve in the class $\ell(\g)$. Thus
$\d(\vp_N,\vp_L)[S^1\x\{ x\}]=\s(\vp_L)_*[S^1\x\{ x\}]\cdot [S^1\x \ell(\g)] = [S^1\x\{ \pt\}]\cdot [S^1\x \ell(\g)] =0$, and
\begin{multline*}
\d(\vp_N,\vp_L)[\{1\}\x\g]= \s(\vp_L)_*[\{1\}\x\g]\cdot [S^1\x \ell(\g)] =\\
[\{1\}\x\g']\cdot [S^1\x \ell(\g)] =\pm (\ell(\g) + \lam(\g) [\bd D^2]) \cdot \ell(\g)
=\pm \lam(\g)
\end{multline*}
Thus $\lam(T_\g)=\text{gcd}(|\d(\vp_N,\vp_L)[\{1\}\x\g]|,|\d(\vp_N,\vp_L)[S^1\x\{ x\}]|) =|\lam(\g)|$.
\end{proof}

 \begin{lem} \label{inf} Given any nontrivial fibered knot $K$ in $S^3$, there is a sequence of embedded loops $\g_n$ contained in a fixed fiber $\Sig$ of $S^3\- K\to S^1$ such that $\lim\limits_{n\to\infty} |\lam(\g_n)|=\infty$.
 \end{lem}
 \begin{proof} If we can find any $c\in H_1(\Sig)$ represented by an embedded loop such that $\lam(c)\ne 0$, then if $e\in H_1(\Sig)$ is represented by a loop and is not a multiple of $c$, $\lam(e+nc)$ is the linking number of $e+nc$ with $e'+nc'$ ($c'$, $e'$ the Lagrangian pushoffs). Thus
 \[ \lam(e+nc)= \lam(e) + n^2\lam(c) + n({\rm{lk}}(c,e')+{\rm{lk}}(e,c')),\]
 whose absolute value clearly goes to $\infty$ as $n\to\infty$. Further, $e+nc$ is represented by an embedded loop for all $n$ for which $e+nc$ is primitive, and this is true for infinitely many $n$. (To see this, identify $H_1(\Sig)$ with $\Z^{2g}$. Then since $c$ and $e$ are independent and primitive, we may make a change of coordinates so that in these new coordinates $c=(1,0,0,\dots,0)$ and $e= (r,s,0,\dots,0)$, for $r,s\in\Z$, $s\ne0$. Thus $e+nc=(n+r,s,0,\dots,0)$. The first coordinate is prime for infinitely many $n$, and at most finitely many of these primes can divide $s$. So these $e+nc$ are primitive.)
 
To find $c$ with $\lam(c)\ne0$, note that ${\rm{lk}}(c,e')$ is the Seifert linking pairing. Since $\DD_K(t)\ne 1$; this pairing is nontrivial. Let $\{ b_i\}$ be a basis for $H_1(\Sig)$. If all $\lam(b_i)=0$, and all $\lam(b_i+b_j)= 0$ then ${\rm{lk}}(b_j,b'_i) = -{\rm{lk}}(b_i,b'_j)$ for all $i\ne j$. This means that the Seifert matrix $V$ corresponding to this basis satisfies $V^T=-V$. However, $\pm1=\DD_K(1)=\det(V^T- V)=\det(2V^T)=2^{2g}\det(V^T)$, a contradiction. 
 \end{proof}

We conclude from Theorem~\ref{LFI} and this lemma:

\begin{thm} \label{sympinv} Let $K$ be any nontrivial fibered knot in $S^3$. Then in the symplectic manifold $X = S^1\x M_K$ there are infinitely many nullhomologous Lagrangian tori which are inequivalent under symplectomorphisms of $X$. \qed
\end{thm}

The constructions of this section are related to Polterovich's `linking class' 
$L\in H^1(T;\Z)$ (see \cite{P}) which is defined for Lagrangian tori $T \C \mathbf{C}^2$, by
$L([a]) = {\rm{lk}}(T, a')$, where $a'$ is a pushoff of a representative $a$ of $[a]\in H_1(T;\Z)$ in the Lagrangian direction. One quickly sees that $L$ is actually defined for a nullhomologous Lagrangian torus in any symplectic $4$--manifold, and $L=\d(\vp_N,\vp_L)$.

The Polterovich linking class is also defined for totally real tori in $\mathbf{C}^2$, and it is shown in \cite{P} that the value of $L$ on totally real tori can be essentially arbitrary, whereas Eliashberg and Polterovich have shown that in $\mathbf{C}^2$ the linking class $L$ vanishes on Lagrangian tori. 

The results of this section may be interpreted as saying that this vanishing phenomenon disappears in symplectic $4$--manifolds more complicated than $\mathbf{C}^2$.

\section{Nullhomologous Lagrangian tori}\label{NHT}

In this section we shall describe examples of collections of $\Ci$\!--inequivalent nullhomologous  Lagrangian tori. The key point is that for our examples, the Lagrangian framing invariant is actually a $\Ci$ invariant.

We begin by describing the symplectic $4$--manifolds which contain the examples.
Let $X$ be a symplectic $4$--manifold with $b_2^+(X)>1$ which contains an embedded symplectic torus $F$ satisfying
\begin{enumerate}
\item[(a)] $F\cdot F=0$
\item[(b)] $F$ contains a loop $\a_{\DD}$, primitive in $\pi_1(F)$, which in $X\- F$ bounds an embedded disk $\DD$ of self-intersection $-1$.
\end{enumerate}
For example, a fiber of a simply connected elliptic surface satisfies this condition. Any torus with a neighborhood symplectically diffeomorphic to a neighborhood of a nodal or cuspidal fiber in an elliptic surface also satisfies the condition, and such tori can be seen to occur in many complex surfaces \cite{rat}.      Let us describe this situation by saying that $X$ contains an embedded symplectic self-intersection $0$ torus {\it with a vanishing cycle}.

Now consider a genus $g$ fibered knot $K$ in $S^3$, and let $\Sig$ be a fiber of the fibration $M_K\to S^1$ and let $m$ be a meridian of $K$.  Let $X$ be a symplectic $4$--manifold with $b_2^+(X)>1$ and with an embedded symplectic self-intersection $0$ torus, $F$, with a vanishing cycle. Fix tubular neighborhoods $N=S^1\x m\x D^2$  of the torus $S^1\x m$ in $S^1\x M_K$ and $N_F=F\x D^2$ of $F$ in $X$,  and consider the result of knot surgery
\[ X_K= X\#_{F=S^1\x m}(S^1\x M_K)\]
where we require the gluing to take the circle $(S^1\x\pt\x\pt)$ in $\bd N$ to $(\a_{\DD}\x\pt)$ in $\bd N_F$. 
Then $X_K$ is a symplectic $4$--manifold with Seiberg--Witten invariant 
 $\sw_{X_K}=\Ds_K(t_F^2)\cdot \sw_X$ where $\Ds_K$ is the symmetrized Alexander  polynomial of $K$ (see \cite{KL4M}). Fix an embedded loop  $\g$ on $\Sig$ whose linking number with the chosen meridian $m$ is $0$, and let $T_\g = S^1\x \g$, a Lagrangian torus with tubular neighborhood $N_{T_\g}= T_\g\x D^2$ in $S^1\x M_K$. Now $M_K$ is a homology $S^1\x S^2$ and $H_1(\Sig)\to H_1(M_K)$ is the $0$--map. Removing the neighborhood $N_m$ of a meridian from $M_K$ does not change $H_1$. ($N_m\cap\Sig =D^2$; so $\bd D^2$ is a meridian to $m$ and it bounds $\Sig\- D^2$.)  Thus
we have $[\g]=0$ in $H_1(M_K\- N_m)$, and hence $T_\g$ is nullhomologous in $X_K$. In fact, since the linking number of $\g$ and $K$ and $m$ is $0$, the loop $\g$ bounds an oriented surface $C\C S^3\- (K\cup m)$. Thus $S^1\x C$ provides a nullhomology of $T_\g$ in $X_K$.
Also note that $b_2^+(X_K\- T_\g)=b_2^+(X)>1$.

\begin{prop} For loops $\g_1$, $\g_2$ in the fiber $\Sig$ of $M_K\to S^1$, if the corresponding nullhomologous tori $T_{\g_1}$ and $T_{\g_2}$ in $X_K$ are symplectically equivalent then $\lam(T_{\g_1}) =\lam(T_{\g_1}) $.
\end{prop}
\begin{proof} Because $S^1\x C$ is a nullhomology of $T_g$, the invariant $\lam(T_\g)$ is calculated exactly as in Theorem~\ref{LFI}; so this proposition follows. \end{proof}

We wish to calculate $\I(X_K,T_\g)$. First fix a basis  for $\bd N_{T_\g}$ which is adapted to the Lagrangian framing of $T_\g$. This basis is $\{ [S^1\x \{ y\}], [\g'], [\bd D^2]\}$ where $\g'$ is a Lagrangian pushoff of $\g$ in $\Sig$ and $y\in \g'$. We begin by studying $X_{K,T_\g}(1,0,0)$, the manifold obtained from $X_K$ by the surgery on $T_\g$ which kills $S^1\x \{ y\}$.
 
\begin {prop}\label{100} $\sw_{X_{K,T_\g}(1,0,0)}=0$. \end{prop}
\begin{proof} Let $\t$ be a path in $\Sig$ from $y$ to the point $x$ at which $m$ intersects $\Sig$. By construction, $S^1\x \{ x\}$ is identified with $\a_\DD\x\pt\in\bd N_F$.
This means that $S^1\x \{ x\}$ is the boundary of a disk $\Delta$ of self-intersection $-1$ in $X\- N_F$. The surgery curve, $S^1\x \{ y\}$, bounds a disk $D$ of self-intersection $0$ in $X_{K,T_\g}(1,0,0)$ (disjoint from $X_K\- T_\g\x D^2$); so the surgered manifold $X_{K,T_\g}(1,0,0)$ contains the sphere $C=\Delta\cup (S^1\x \t)\cup D$ of self-intersection $-1$.  

The rim torus $R=m\x \bd D^2\C\bd N_F$ intersects the sphere $C$ in a single positive intersection point, but this is impossible if $\sw_{X_{K,T_\g}(1,0,0)}\ne 0$. For, if $\sw_{X_{K,T_\g}(1,0,0)}\ne0$, then blowing down $C$, we obtain a $4$--manifold $Z$ (with $b_2^+>1$) which contains a torus $R'$ of self-intersection $+1$, and the  Seiberg--Witten invariant of $Z$ is nontrivial.  However, the adjunction inequality states that for any basic class $\b$ of $Z$ we have $0\ge 1 +|\b\cdot R'|$, an obvious contradiction.  \end{proof}
  
Before proceeding further, note that since $T_\g$ is nullhomologous in $X_K$,
\[ j_*\co H_2(X_K)\to H_2(X_K,T_\g)\]
is an injection. Thus we may identify $\swb_{(X_K,T_\g)}=j_*(\sw_{X_K})$ with $\sw_{X_K}$.  We shall make use of an important result due to Meng and Taubes concerning the Seiberg--Witten invariant of a closed $3$--manifold $M$ \cite{MT}:
\begin{equation}\label{MTF}
\begin{cases}
\sw_M=\Ds_M(t^2)\cdot (t-t^{-1})^{-2}, \ \ b_1(M)=1\\
\sw_M=\Ds_M, \ \  b_1(M)>1
\end{cases}\end{equation}
where $\Ds_M$ is the symmetrized Alexander polynomial of $M$, and if $b_1(M)=1$ then $t\in \Z{H_1(M;\R)}$ corresponds to the generator of $H_1(M,\R)$. 

Since $X_{K,T_\g}(0,0,1)$ is the result of the surgery which kills $\bd D^2$, it is $X_K$ again, and we know that $\sw_{X_K} = \Ds_K(t_F^2)\cdot\sw_X$. 
This also means that $\swb_{(X_K,T_\g)} = \Ds_K(t_F^2)\cdot\sw_X$. Thus to calculate $\I(X_K,T_\g)$, it remains only to calculate the Seiberg--Witten invariant of $X_{K,T_\g}(0,1,0)$, the manifold obtained by the surgery on $T_\g$ which makes $\g'$ bound a disk. 

Let $M_K(\g)$ denote the result of surgery on $\g$ in $M_K$ with the Lagrangian framing. In terms of the usual nullhomologous framing, this is the result of surgery on the link $K\cup \g$ in $S^3$ with framings $0$ on $K$ and $\lam(\g)$ on $\g$. 
In case $\lam(\g)\ne 0$, we have $b_1(M_K(\g))=1$ and if $\lam(\g)=0$ then $b_1(M_K(\g))=2$. In this case, the extra generator of $H_1(M_K(\g);\R)$ is given by a meridian to $\g$ in $S^3$.
Accordingly, the Seiberg--Witten invariant of $M_K(\g)$ (equivalently, the Seiberg--Witten invariant of $S^1\x M_K(\g)$) is given by
\begin{equation} \label{swmkg}   \sw_{M_K(\g)} = 
\begin{cases}
\Ds_{M_K(\g)}(t^2)\cdot (t-t^{-1})^{-2}, \ \ \lam(\g)\ne 0 \\
\Ds_{M_K(\g)}(t^2,s^2), \ \ \lam(\g)=0
\end{cases}
\end{equation}
where $t$ corresponds to the meridian of $K$ and $s$ to the meridian of $\g$.

\begin{prop} \label{order} Suppose that $\lam(\g)\ne 0$, then $|\DD_{M_K(\g)}(1)| = |\lam(\g)|$.
\end{prop}
\begin{proof} We have $H_1(M_K(\g)) = \Z\oplus\Z_{|\lam(\g)|}$. It is a well-known fact \cite{Turaev} that for $3$--manifolds with $b_1=1$, the sum of the coefficients of the Alexander polynomial is, up to sign, the order of the torsion of $H_1$.
\end{proof}

Let $Z$ denote the $3$--manifold obtained from $M_K$ by doing $+1$ surgery on $T_\g$ with respect to the Lagrangian framing. This is the surgery that kills the class $[\g']+[\bd D^2]$ on the boundary of a tubular neighborhood $\g \x D^2$ of $\g$. Then $H_1(Z)=\Z \oplus\Z_b$ where $b=|\lam(\g)+1|$.
Since $M_K$ is fibered over the circle, the manifold $Z$ is also fibered over the circle with the same fiber $\Sig$.
(This is true for any $(1/p)$--Lagrangian-framed surgery. The effect of such a surgery on the monodromy is to compose it with the $p$th power of a Dehn twist about $\g$. See \cite{S, ADK}.)

If $\lam(\g)\ne -1$, then $b_1(Z)=1$, and its symmetrized Alexander polynomial $\Ds_Z(t)$ is a function of one variable. If  $\lam(\g)=-1$, we have $H_1(Z)=\Z\oplus\Z$.  In this case, the Alexander polynomial of $Z$ is a $2$--variable polynomial $\Ds_Z(t,s)$ where $s$ corresponds to the meridian of $\g$. Let $\bar{\DD}^{sym}_Z(t)=\Ds_Z(t,1)$. 
 
Write $\Ds_K(t)=a_0+ (t^{g} +t^{-g})+\sum\limits_{i=1}^{g-1} a_i(t^i+t^{-i})$. (This is equal to $\Ds_{M_K}(t)$.)

\begin{lem} \label{fibered}  The symmetrized Alexander polynomial of $Z$ is given by 
\begin{gather*}\Ds_{Z}(t)=b_0+ (t^{g} +t^{-g}) +\sum\limits_{i=1}^{g-1} b_i(t^i+t^{-i}),\ \ 
\lam(\g)\ne -1\\
\bar{\DD}^{sym}_Z (t)= \left(b_0+ (t^{g} +t^{-g}) +\sum\limits_{i=1}^{g-1} b_i(t^i+t^{-i})\right)\cdot
(t^{1/2}-t^{-1/2})^{-2}, \ \ \lam(\g)=-1\end{gather*}
for some choice of coefficients $b_i$.
\end{lem}
\begin{proof}
For a $3$--manifold with $b_1=1$ which is fibered over the circle, the Alexander polynomial is the characteristic polynomial of the (homology) monodromy. (Compare \cite[VII.5.d]{Turaev}.) This is a monic symmetric polynomial of degree $2g$, as claimed.

In case $\lam(\g)=-1$, one can either apply the theorem of Turaev {\it{op.cit.}} or apply the work of Hutchings and Lee. According to \cite{HL} together with Mark \cite{M}, after appropriate symmetrization, the zeta invariant of the monodromy, namely the characteristic (Laurent) polynomial of the homology monodromy times the term $(t^{1/2}-t^{-1/2})^{-2}$ is equal to a (Laurent) polynomial in $t$, whose coefficient of $t^n$ is the sum over $m$ of the coefficients of all terms of $\Ds_Z(t,s)$ of the form $a_{n,m}t^ns^m$. In other words, 
$\bar{\DD}^{sym}_Z (t)=\Ds_Z(t,1)$ is this Laurent polynomial. This proves the second statement of the lemma.
\end{proof}

\begin{lem} \label{DD} The Seiberg--Witten invariant of $X_{K,T_\g}(0,1,0)$ is
 \[\swb_{(X_{K,T_\g}(0,1,0),T_\g)} =
 \left(\Ds_{Z}(t_F^2)-\Ds_K(t_F^2)\right)\cdot\sw_X\]  if $\lam(\g)\ne -1$,
 and if $\lam(\g)=-1${\em :}
   \[ \swb_{(X_{K,T_\g}(0,1,0),T_\g)} =\left(\bar{\DD}^{sym}_{Z}(t_F^2)\cdot (t_F-t_F^{-1})^2-\Ds_K(t_F^2)\right)\cdot\sw_X, \]
 \end{lem}
\begin{proof} The result of $(+1)$--Lagrangian-framed surgery on $T_\g$ in $X_K$ is the fiber sum $X\#_{F=S^1\x m}(S^1\x Z)$. If $\lam(\g)\ne -1$, it follows from \eqref{MTF} and the usual gluing formulas that 
this manifold has Seiberg--Witten invariant equal to $\Ds_{Z}(t_F^2)\cdot\sw_X$.
Applying the surgery formula, 
 \[ \swb_{(X\#_{F=S^1\x m}(S^1\x Z),T_\g)} = \swb_{(X_{K,T_\g}(0,1,0),T_\g)} +\swb_{(X_K,T_\g)}\]
 or 
\[
\Ds_{Z}(t_F^2)\cdot\sw_X=\swb_{(X_{K,T_\g}(0,1,0),T_\g)}+\Ds_K(t_F^2)\cdot\sw_X
 \]
 and the lemma follows.
 
 If $\lam(\g)=-1$, then 
 \begin{multline*} \swb_{(X\#_{F=S^1\x m}(S^1\x Z),T_\g)}=\sw_X\cdot\swb_{(S^1\x Z,T_\g)} \cdot (t_F -t_F^{-1})^2\\ =
\sw_X\cdot\bar{\DD}^{sym}_Z(t_F^2)\cdot (t_F -t_F^{-1})^2 \end{multline*} 
and the result follows as above.
 \end{proof}

\begin{thm}\label{invt} Let $X$ be a symplectic $4$--manifold with $b^+>1$ containing an embedded self-intersection $0$ torus $F$ with a vanishing cycle. Let $K$ be a nontrivial fibered knot, and let $\g$ be an embedded loop on a fiber of $S^3\- K\to S^1$. Then the Lagrangian framing invariant $\lam(T_\g)$ is an orientation-preserving diffeomorphism invariant of the pair $(X_K, T_\g)$.
 \end{thm}
 \begin{proof} Using the notation above and Lemmas~\ref{fibered} and \ref{DD},
 \begin{equation}\label{SWW} \swb_{(X_{K,T_\g}(0,1,0),T_\g)}= \left( (b_0-a_0)+\sum\limits_{i=1}^{g-1} (b_i-a_i) (t_F^{2i}+t_F^{-2i})\right)\cdot\sw_X 
 \end{equation}
 It follows from (\ref{surgery formula 2}) that
 \begin{multline}\label{0pq}
\swb_{(X_{K,T_\g(0,p,q)},T_\g)}= p\left( (b_0-a_0)+\sum\limits_{i=1}^{g-1} (b_i-a_i)(t_F^{2i}+t_F^{-2i})\right)\cdot\sw_X\\
  +q \left(  a_0+ (t_F^{2g} +t_F^{-2g})+ \sum\limits_{i=1}^{g-1} a_i(t_F^{2i}+t_F^{-2i}) \right)\cdot\sw_X
 \end{multline}
 Since $\Ds_K(1)=\pm1$, we have $2+a_0+2\sum\limits_1^{g-1}a_i=\ve=\pm1$.
Furthermore, 
\begin{equation}\label{anotherSW}\swb_{(X_{K,T_\g}(0,1,0),T_\g)}= \begin{cases} 
\sw_X\cdot \Ds_{M_K(\g)}(t_F^2), \ \ \lam(\g)\ne 0,\\
\sw_X\cdot \bar{\DD}^{sym}_{M_K(\g)}(t_F^2)\cdot (t_F-t_F^{-1})^2, \ \ \lam(\g)=0.
\end{cases}
\end{equation}
(See \eqref{swmkg}.) Thus if $\lam(\g)\ne 0$,
\[ \Ds_{M_K(\g)}(t^2)= (b_0-a_0)+\sum\limits_{i=1}^{g-1} 
(b_i-a_i) (t^{2i}+t^{-2i})\]
and by Proposition \ref{order}, 
\begin{equation}\label{finally}(b_0-a_0)+2\sum\limits_1^{g-1}(b_i-a_i)=\pm\lam(\g)
 =\d\,\lam(\g)\end{equation}
If $\lam(\g)=0$, it follows from equations~\eqref{SWW} and \eqref{anotherSW} that
 \[ \bar{\DD}^{sym}_{M_K(\g)}(t^2)\cdot (t-t^{-1})^2= (b_0-a_0)+\sum\limits_{i=1}^{g-1} 
(b_i-a_i) (t^{2i}+t^{-2i})\]
so $(b_0-a_0)+2\sum\limits_1^{g-1}(b_i-a_i)=0$ in this case, and we see that \eqref{finally} holds in general.
 
 Let $\s(p,q)$ be the sum of all coefficients of $\swb_{(X_{K,T_\g}(0,p,q),T_\g)}/\sw_X$ from terms of degree not equal to $\pm 2g$. Then it follows from \eqref{0pq} that $\s(p,q)= p\d\lam(\g)+q(\ve-2)$. Let $\t(p,q)$ be the coefficient of $t_F^{2g}$ in  
 $\swb_{X_{T_\g}(0,p,q)}/\sw_X$; so $\t(p,q)=q$.
 
 We have seen that \[\I(X_K,T_\g)=\{\swb_{(X_{K,T_\g}(a,p,q),T_\g)}|a,p,q\in\Z\}
 =\{\swb_{(X_{K,T_\g}(0,p,q),T_\g)}|p,q\in\Z\}\]
(the last equality by Proposition \ref{100}) is an orientation-preserving diffeomorphism invariant of the pair $(X_K,T_\g)$. From $\I(X_K,T_\g)$ we can extract the invariant 
\[ \gcd \{|\s(p,q)+(2-\ve)\t(p,q)|\,\big\vert p,q\in\Z\}=\gcd\{|p\d\lam(\g)|\,\big\vert p\in\Z\} =|\lam(\g)|=\lam(T_\g). \]   
\end{proof}

 \begin{thm} \label{NLT}  Let $X$ be a symplectic $4$--manifold with $b^+>1$ containing an embedded self-intersection $0$ torus $F$ with a vanishing cycle, and let $K$ be a nontrivial fibered knot. Then in $X_K$ there is an infinite sequence of pairwise inequivalent nullhomologous Lagrangian tori $T_{\g_n}$. 
\end{thm}
 \begin{proof} Choose a sequence of loops $\g_n$ as in the statement of  Lemma~\ref{inf}, then it is clear that the elements $\g_n$ give inequivalent $T_{\g_n}$.
 \end{proof}

\section{Circle sums and essential Lagrangian tori}

We next discuss a technique for utilizing families of inequivalent nullhomologous Lagrangian tori to build families of inequivalent essential homologous Lagrangian tori. For a fibered knot $K$ in $S^3$, let $m_1$, $m_2$ be meridians, and consider a meridian $\mu_1$ of $m_1$ in $M_K$ as shown in Figure~3.

Let $X_i$, $i=1,2$ be symplectic $4$--manifolds, containing embedded symplectic tori $F_i$ of self-intersection $0$, and suppose that $F_1$ has a vanishing cycle.
Denote by $X$ the fiber sum $X= X_1\#_{F_1=F_2}X_2$ . Then
\[X_K=(X_1\#_{F_1=F_2}X_2)_K=X_1\#_{F_1=S^1\x m_1}S^1\x M_K \#_{S^1\x m_2=F_2}X_2\] 
As in the previous section, we insist that the gluing map from the boundary of a tubular neighborhood of $S^1\x m_1$ to the boundary of a tubular neighborhood of $F_1$ should take $S^1\x\pt\x\pt$ to the loop $\a_\DD\x\pt$ representing the vanishing cycle.
Because $\tm =S^1\x \mu_1$ is a rim torus to $F_1$ 
it follows easily that $\tm$ is a Lagrangian torus which represents an essential, in fact primitive, class in $H_2(X_K)$. 

 \centerline{\small\unitlength 1cm
\begin{picture}(6,5.7)
\put (3,3){\oval(3,4)[r]}
\put (2,5){\line(3,-4){1}}
\put (3,3.66){\line(-3,-4){.4}}
\put (3,5){\line(-3,-4){.4}}
\put (2,3.66){\line(3,4){.4}}
\put (2,3.66){\line(3,-4){1}}
\put (2,2.33){\line(3,4){.4}}
\put (2,2.33){\line(3,-4){1}}
\put (3,2.33){\line(-3,-4){.4}}
\put (2,1){\line(3,4){.4}}
\put (2.4,3.66){\oval(1.5,1.4)[l]}
\put (2.6,3.66){\oval(1.5,1.4)[r]}
\put (1.5,4.25){\small{$\g_1$}}
\curve(2,5,1,5,0.6,4.75,0.5,4)
\curve(.4,4.3, .35,4.29, .3,4.27, .25,4.22, .225,4.15, .224,4.05, .275,3.95, .35,3.92, .4,3.915, .52,3.915, .57,3.92, .645,3.95, .696,4.05, .697,4.15)
\curve(.62,4.08,.58,4.13,.57,4.17,.58,4.21,.62,4.26,.75,4.26,.79,4.21,.80,4.17,.79,4.13,.75,4.08)
\curve(.535,4.3,.58,4.28)
\curve(.5,3.8,.5,2)
\curve(.4,2.25, .35,2.24, .3,2.22, .25,2.17, .225,2.1, .224,2, .275,1.9, .35,1.87, .4,1.865, .52,1.865, .57,1.87, .645, 1.9, .696,2, .695,2.1, .67,2.17,.62,2.22, .57,2.25 )
\curve(2,1,1,1,0.6,1.25,.5,1.75)
\put(-.25,3.7){\small{$m_1$}}
\put(-.25,1.7){\small{$m_2$}}
\put(.75,4.4){\small{$\mu_1$}}
\put(2.5,4.75){\small{$0$}}
\put (1.85,.5){Figure 3}
\end{picture}}

\begin{prop} $\I(X_K,\tm)=\{q\,\swb_{(X_K,T_{\mu})}|\,q\in\Z\}$. \end{prop}
\begin{proof}  If we choose a tubular neighborhood $\tm\x D^2$ with the Lagrangian framing, then we may use the basis 
$\{[S^1\x \{ y\}], [\mu_1'],[\bd D^2]\}$ for $H_1(\bd(\tm\x D^2))$, and as in Proposition~\ref{100},  $\sw_{X_{K,\tm}(1,0,0)}=0$. Of course, $\sw_{X_{K,\tm}(0,0,1)}=\sw_{X_K}$.  It remains to calculate $\sw_{X_{K,\tm}(0,1,0)}$.
This can be done via Kirby calculus. Let $Y$ denote the result of surgery on $\mu_1$ in $M_K$ with respect to the Lagrangian framing. This is shown in  Figure~4.

Slide the $0$--surgered handle corresponding to $K$ over the $0$--surgered $\mu_1$, to get Figure~5. Thus $Y\cong M_K\# (S^1\x S^2)$. Hence
\begin{multline*}X_{K,\tm}(0,1,0)=(X_1\#_{F_1=S^1\x m_1} T^2\x S^2)\- (D^3\x S^1)\\ \bigcup_{S^2\x S^1}\, (X_2\#_{F_2=S^1\x m_2} S^1\x M_K)\- (D^3\x S^1)
\end{multline*}
Each of the two sides of the union has $b^+\ge1$, and $S^2\x S^1$ admits a metric of positive scalar curvature. This means that $\sw_{X_{K,\tm}(0,1,0)}=0$, and the proposition follows.
\end{proof}
 
\centerline{\small\unitlength 1cm
\begin{picture}(12,5)
\put (3,3){\oval(3,4)[r]}
\put (2,5){\line(3,-4){1}}
\put (3,3.66){\line(-3,-4){.4}}
\put (3,5){\line(-3,-4){.4}}
\put (2,3.66){\line(3,4){.4}}
\put (2,3.66){\line(3,-4){1}}
\put (2,2.33){\line(3,4){.4}}
\put (2,2.33){\line(3,-4){1}}
\put (3,2.33){\line(-3,-4){.4}}
\put (2,1){\line(3,4){.4}}
\put (2.4,3.66){\oval(1.5,1.4)[l]}
\put (2.6,3.66){\oval(1.5,1.4)[r]}
\put (1.5,4.25){\small{$\g_n$}}
\curve(2,5,1,5,0.6,4.75,0.5,4)
\curve(.4,4.3, .35,4.29, .3,4.27, .25,4.22, .225,4.15, .224,4.05, .275,3.95, .35,3.92, .4,3.915, .52,3.915, .57,3.92, .645,3.95, .696,4.05, .697,4.15)
\curve(.62,4.08,.58,4.13,.57,4.17,.58,4.21,.62,4.26,.75,4.26,.79,4.21,.80,4.17,.79,4.13,.75,4.08)
\curve(.535,4.3,.58,4.28)
\curve(.5,3.8,.5,2)
\curve(.4,2.25, .35,2.24, .3,2.22, .25,2.17, .225,2.1, .224,2, .275,1.9, .35,1.87, .4,1.865, .52,1.865, .57,1.87, .645, 1.9, .696,2, .695,2.1, .67,2.17,.62,2.22, .57,2.25 )
\curve(2,1,1,1,0.6,1.25,.5,1.75)
\put(-.25,3.7){\small{$m_1$}}
\put(-.25,1.7){\small{$m_2$}}
\put(.775,4.3){\small{$0$}}
\put(3.5, 4.65){\small{$0$}}
\put (1.85,.5){Figure 4}
\curve(7.58,4.28,7.4,4.3, 7.35,4.29, 7.3,4.27, 7.25,4.22, 7.225,4.15, 7.224,4.05, 7.275,3.95, 7.35,3.92, 7.4,3.915, 7.52,3.915, 7.57,3.92, 7.645,3.95, 7.696,4.05, 7.697,4.15)
\curve(7.62,4.08,7.58,4.13,7.57,4.17,7.58,4.21,7.62,4.26,7.75,4.26,7.79,4.21,
7.8,4.17,7.79,4.13,7.75,4.08)
\put(6.75,3.7){\small{$m_1$}}
\put(7.775,4.3){\small{$0$}}
\put (11,3){\oval(3,4)[r]}
\put (10,5){\line(3,-4){1}}
\put (11,3.66){\line(-3,-4){.4}}
\put (11,5){\line(-3,-4){.4}}
\put (10,3.66){\line(3,4){.4}}
\put (10,3.66){\line(3,-4){1}}
\put (10,2.33){\line(3,4){.4}}
\put (10,2.33){\line(3,-4){1}}
\put (11,2.33){\line(-3,-4){.4}}
\put (10,1){\line(3,4){.4}}
\curve(10,5,9,5,8.6,4.75,8.5,4)
\curve(.4,4.3, .35,4.29, .3,4.27, .25,4.22, .225,4.15, .224,4.05, .275,3.95, .35,3.92, .4,3.915, .52,3.915, .57,3.92, .645,3.95, .696,4.05, .697,4.15)
\curve(8.5,4,8.5,2)
\curve(8.4,2.25, 8.35,2.24, 8.3,2.22,8.25,2.17, 8.225,2.1, 8.224,2, 8.275,1.9, 8.35,1.87, 8.4,1.865, 8.52,1.865, 8.57,1.87,8.645, 1.9, 8.696,2, 8.695,2.1, 8.67,2.17,8.62,2.22, 8.57,2.25 )
\curve(10,1,9,1,8.6,1.25,8.5,1.75)
\put(7.75,1.7){\small{$m_2$}}
\put(11.5, 4.65){\small{$0$}}
\put (9,.5){Figure 5}
\end{picture}}
 
Let $\{ T_{\g_n} \}$ be the family of nullhomologous Lagrangian tori in $X_K$ given by Theorem~\ref{NLT} (thinking of $m=m_1$). For any $T_{\g_n}$ in the family, we can form the circle sum of 
$T_{\g_n}$ with $\tm$. This is done by fixing a path in the fiber $\Sig$ running from a point of $\g_n$ to a point in $\mu_1$ and taking the connected sum $\g_n'$ of $\g_n$ with $\mu_1$ along this path. (See Figures~4 and 6.) The resulting torus $T_{\g_n'}=S^1\x\g_n'$ in $X$ is Lagrangian and homologous to $T_{\mu}$. Note that $\lam(\g_n')=\lam(\g_n)$; in fact, $\g'_n$ and $\g_n$ are isotopic in $M_K$.

\centerline{\small\unitlength 1cm
\begin{picture}(6,5.5)
\put (3,3){\oval(3,4)[r]}
\put (2,5){\line(3,-4){1}}
\put (3,3.66){\line(-3,-4){.4}}
\put (3,5){\line(-3,-4){.4}}
\put (2,3.66){\line(3,4){.4}}
\put (2,3.66){\line(3,-4){1}}
\put (2,2.33){\line(3,4){.4}}
\put (2,2.33){\line(3,-4){1}}
\put (3,2.33){\line(-3,-4){.4}}
\put (2,1){\line(3,4){.4}}
\put (2.6,3.66){\oval(1.5,1.4)[r]}
\put (1.5,4.25){\small{$\g_n'$}}
\curve(2,5,1,5,0.6,4.75,0.5,4)
\curve(.4,4.3, .35,4.29, .3,4.27, .25,4.22, .225,4.15, .224,4.05, .275,3.95, .35,3.92, .4,3.915, .52,3.915, .57,3.92, .645,3.95, .696,4.05, .697,4.15)
\curve(.62,4.08,.58,4.13,.57,4.17,.58,4.21,.62,4.26,.75,4.26,.79,4.21,.80,4.19, .9, 4.12,1.62,4,1.8,4.07,2,4.225,2.2,4.32,2.4,4.35)
\curve(.725,4.05,.86,4.02,.95,4,1.1, 3.97, 1.5,3.9,1.6, 3.87,1.65,3.8,1.76,3.6,1.8,3.5,1.9,3.34,2,3.23,2.05,3.17,2.13,3.1,2.27,3,2.35,2.96,2.34,2.97)
\curve(.535,4.3,.58,4.28)
\curve(.5,3.8,.5,2)
\curve(.4,2.25, .35,2.24, .3,2.22, .25,2.17, .225,2.1, .224,2, .275,1.9, .35,1.87, .4,1.865, .52,1.865, .57,1.87, .645, 1.9, .696,2, .695,2.1, .67,2.17,.62,2.22, .57,2.25 )
\curve(2,1,1,1,0.6,1.25,.5,1.75)
\put(-.25,3.7){\small{$m_1$}}
\put(-.25,1.7){\small{$m_2$}}
\put (1.85,.5){Figure 6}
\end{picture}}

We wish to calculate $\I(X_K,T_{\g'_n})$. Using, as in the previous section, the basis adapted to the Lagrangian framing of $T_{\g'_n}$, because $\sw_{X_{K,T_{\g'_n}}(1,0,0)}=0$ and $\sw_{X_{K,T_{\g'_n}}(0,0,1)}=\sw_{X_K}$,
we only need to calculate $\swb_{(X_{K,T_{\g'_n}}(0,1,0),T_{\g'})}$.

Consider $M_K(\g_n')$ the $3$--manifold obtained by surgery on $\g'_n$ in $M_K$ using the Lagrangian framing. Since $\g_n'$ and $\g_n$ are isotopic in $M_K$, the manifolds $M_K(\g_n')$ and $M_K(\g_n)$ are diffeomorphic. We may as well assume that $\lam(\g)\ne 0$.
We have $H_1(M_K(\g_n'))=\Z\oplus\Z_{|\lam(\g_n)|}$ where the infinite cyclic summand is generated by the class of $m_2$ and the finite cyclic summand by the class of $\nu$, a meridian to $\g_n$. Note that $[m_1]=[m_2]+[\nu]$. 

To obtain $M_K(\g_n')$, one does surgery which kills the curve $\g_n'+\lam(\g_n)\,\nu$. Thus $[\g_n'] = -\lam(\g_n)\,[\nu]$ in 
$H_1(M_K(\g_n')\- (m_1\cup m_2))$.
Furthermore, in the manifold with boundary $M_K(\g_n')\- (m_1\cup m_2)$, the core $C$ of the surgery solid torus is homologous to $(1/\lam(\g_n))\,\g_n'$. We are interested in
\[ X_{K,T_{\g'_n}}(0,1,0) = X_1\#_{F_1=S^1\x m_1} S^1\x M_K(\g_n') \#_{S^1\x m_2=F_2}X_2\]
In $H_2(X_{K,T_{\g'_n}}(0,1,0))$ we have 
\[  [F_1] = [F_2] + [S^1\x \nu] = [F_2] - (1/\lam(\g_n)) [S^1\x \g_n'] 
= [F_2]- [S^1\x C] \]
Let $T= S^1\x C$. We need to calculate $\swb_{(X_{K,T_{\g'_n}}(0,1,0),T_{\g_n'})}$ which is an element of 
\[ \Z H_2(X_{K,T_{\g'_n}}(0,1,0)\- T\x D^2,\bd) =
\Z H_2(X_{K,T_{\g'_n}}(0,1,0),T)\]
Thus we may assume that 
$[F_1] = [F_2]=[F]$, say, for the purpose of this calculation.
Now precisely the same proof as that of Theorem~\ref{invt} (replacing $X$ by $X_1\#_{F_1=F_2}X_2$) gives:

\begin{thm} \label{ess} Let $X_i$, $i=1,2$ be symplectic $4$--manifolds containing embedded symplectic tori $F_i$ of self-intersection $0$. Suppose also that $F_1$ has a vanishing cycle. Set $X=X_1\#_{F_1=F_2}X_2$. Let $K$ be a nontrivial fibered knot and let $\g$ be an embedded loop on a fiber $\Sig$ of $S^3\- K\to S^1$.  Let $m_1$ and $m_2$ be meridians to $K$ which do not link $\g$,  $\mu$ a meridian to $m_1$ which lies on $\Sig$, and $\g'$ the connected sum of $\g$ and $\mu$ in $\Sig$. Then $T_{\g'}=S^1\x\g'$ is a Lagrangian torus in $X_K$ and represents a primitive homology class.

 If $\eta$ is another loop on $\Sig$ which has linking number $0$ with $m_1$ and $m_2$ and $\eta'$ is  the connected sum of $\eta$ and $\mu$ in $\Sig$ with corresponding Lagrangian torus $T_{\mu'}=S^1\x \mu'$, then $T_{\eta'}$ is homologous to $T_{\g'}$ and if  $T_{\g'}$ and $T_{\eta'}$ are equivalent in $X_K$, it follows that $\lam(T_\eta)=\lam(T_\g)$.  \qed
 \end{thm}
 
 This completes the proof of Theorem~\ref{T1}. Since the fiber $F$ of $E(1)$ has a vanishing cycle and since $(E(1)\#_FX)_K\cong X_K\#_FE(1)$ we have, for example, the following corollary:
 
\begin{cor}\label{ex} Let $X$ be any symplectic $4$--manifold containing an embedded symplectic torus of self-intersection $0$. Let $K$ be any nontrivial fibered knot. Then in $X_K\#_FE(1)$ there is an infinite sequence of essential Lagrangian tori $T_{\g'_n}$ which are pairwise homologous but no two of which are equivalent. \qed
\end{cor}

\section{Symplectically inequivalent Lagrangian tori which are not distinguished by relative Seiberg--Witten invariants}\label{SELT}

These examples live in simple versions of the manifolds constructed in \cite{1bc}. It would be very easy to give much larger classes of examples, but we shall content ourselves with those below.  

Let $X$ be the K3--surface with elliptic fiber $F$, and let $K$ be a nontrivial fibered knot. Let $\g_1$, $\g_2$  be embedded loops on a fiber of $S^3\- K\to S^1$ with different $\lam(\g_i)>0$, and consider the nullhomologous Lagrangian tori $T_{\g_i}\C X_K$. Theorem~\ref{invt} implies that these tori are $\Ci$\!--inequivalent. Suppose that the knot $K$ has genus $g$. In the construction of
$X_K$ we have replaced a 2--disk in a section of the elliptic fibration of K3 
with a punctured surface of genus $g$, a fiber of $S^3\- K\to S^1$. Thus $X_K$ contains a symplectic genus $g$ surface $S$ of self-intersection $-2$ and $[S]\cdot [F] = 1$. Consider another smooth fiber $F'$ of the elliptic fibration of $(X\setminus N(T))\subset X_K$. Then $F'+S$ is a singular surface with one double point, which can be smoothed to give a symplectically embedded surface $S'$ of genus $g+1$ representing the homology class $[S']=[F]+[S]$. Thus $[S']^2=0$ and $[S']\cdot [F]=1$. 
 
Next, let $K'$ denote the trefoil knot in $S^3$. Since $K'$ is a fibered genus $1$ knot, the 4--manifold $S^1\x M_{K'}$ is a smooth $T^2$--fiber bundle over $T^2$. We obtain a symplectic manifold $Y$ by forming the fiber sum of $g+1$ copies of $S^1\x M_{K'}$ along the tori $T$ where $m_i$ is a meridian of $K'$ in the $i$th copy of $M_{K'}$. There is a fiber bundle 
\[ \begin{array}{ccc}
T^2 \longrightarrow & Y&\negthickspace =\ S^1\x M_{K'}\#_{T^2}\cdots\#_{T^2}S^1\x M_{K'}\\
&\Big\downarrow&\\
& C_0&
\end{array}\]
where $C_0$ is a genus $g+1$ surface.  There is a symplectic section $C\subset Y$ given by the connected sum of the individual sections $S^1\x{m_i}$. 
Now form the symplectic manifold $Z_K=X_K\#_{S'=C}Y$. We perform the fiber sum so that $Z_K$ is spin \cite{Gompf}. In \cite{1bc} it is shown that $Z_K$ is simply connected and that $\sw_{Z_K}=t_\k+(-1)^g t_\k^{-1}$ where $\k$ is the canonical class of $Z_K$.

Since $S'\C X_K$ can be chosen disjoint from a nullhomology of $T_{\g_i}$, we still have $\lam(T_{\g_i})=\lam(\g_i)$ for the Lagrangian framing invariants of the nullhomologous Lagrangian tori $T_{\g_i}\C Z_K$. Thus $T_{\g_1}$ and $T_{\g_2}$ are symplectically inequivalent in $Z_K$.  

We next compute $I(Z_K,T_{\g_i})$, using the same basis that we used in \S\ref{NHT}.
As in Proposition~\ref{100}, $\sw_{Z_{K,T_\g}(1,0,0)}=0$. The point is that 
\[ Z_{K,T_\g}(1,0,0)=X_{K,T_\g}(1,0,0)\#_{S'=C}Y\] and
the sphere of self-intersection $-1$ and its dual torus of square $0$ found in Proposition~\ref{100} both live in the complement of the surface $S'$.

Next consider $Z_{K,T_\g}(0,1,0)$. For any manifold $\bX$ obtained from $X$ by surgery on $T_\g$, $H_2(\bX)=H_2(X)$ if the surgery curve is homologically nontrivial in 
$X\- T_\g$ and $H_2(\bX)=H_2(X)\oplus U$ where $U$ has rank $2$ if the surgery curve is homologically trivial in $X\- T\g$. Since $\lam(\g_i)>0$, the surgery curve $\g'$ for the surgery giving $Z_{K,T_\g}(0,1,0)$ is not nullhomologous. This means that 
for 
\[ Z_{K,T_\g}(0,1,0)=X_{K,T_\g}(0,1,0)\#_{S'=C}Y\]
$H_2(Z_{K,T_\g}(0,1,0)) = H_2(Z_K)$, and the arguments of \cite{1bc} using the adjunction inequality and \cite{MST} to show that $\sw_{Z_K}=t_\k+(-1)^g t_\k^{-1}$, will again show that $\sw_{Z_{K,T_\g}(0,1,0)}=t_\k+(-1)^g t_\k^{-1}$.
Hence $I(Z_K,T_{\g_i})=\{p(t_\k+(-1)^g t_\k^{-1})+q(t_\k+(-1)^g t_\k^{-1})|p,q\in \Z\}=\{r(t_\k+(-1)^g t_\k^{-1})|r\in \Z\}$, independent of $\g_i$; so relative Seiberg--Witten invariants don't detect whether or not the $T_{\g_i}$ are $\Ci$\!--equivalent.
It would be extremely interesting if they were not.

\section{Discussion}

As we have already mentioned in \S \ref{Intro}, the first examples like those of Corollary~\ref{ex} were recently discovered by S. Vidussi \cite{V}. The examples of \cite{V} live in $E(2)_K$ and are of the type described in Theorem~\ref{ess}, thus they can be distinguished by using the Lagrangian framing invariant. Again view $E(2)_K$ as a double fiber sum, $E(1)\#_{F=S^1\x m_1}S^1\x M_K \#_{S^1\x m_2=F}E(1)$.
Vidussi points out that for a torus $T_{\b} = S^1\x\b$, $\b$ a loop in $\Sig$ whose linking number with $m_2$ is $0$, the homology class of $T_{\b}$ in $E(2)_K$ is determined by the linking number of  $\b$ with $m_1$. We have restricted ourselves to the case where the linking number is $1$, but this is completely unnecessary, and the Lagrangian framing invariant gives an invariant for the general situation.

One need not restrict to genus one Lagrangian submanifolds in order to take advantage of the technique of Lagrangian circle sums with nullhomologous Lagrangian tori.  However the authors have not yet been able to find invariants for higher genus Lagrangian surfaces which are as simple to calculate as those in this paper.

\rk{Acknowledgements}The authors gratefully acknowledge support from the National Science Foundation. The first author was partially supported NSF Grants DMS0072212 and DMS0305818, and the second author by NSF Grant DMS0204041.

\end{document}